%%%%%%%%%%%%%%%%%%%%%%%%%%%%%%%%%%%%%%%%%%%%%%%%%%%%%%%%%%%%%%%%%%%%%%%%%%%
%%% Title:  Characteristic classes for $G$-structures
%%% Author: D\. V\. Alekseevsky, Peter W. Michor
%%% Remark: AmSTeX, 7 pages
%%% Series: Differential Geometry and its Applications 3 (1993) 323--329
%%%%%%%%%%%%%%%%%%%%%%%%%%%%%%%%%%%%%%%%%%%%%%%%%%%%%%%%%%%%%%%%%%%%%%%%%%%
% TeX-NMB program applied to the file from 1991.11.17;13:38 on 1991.11.17; 13:39
% TeX-NMB program applied to the file from 1991.11.1;9:9 on 1991.11.1; 9:9
% TeX-NMB program applied to the file from 1991.6.15;12:19 on 1991.6.15; 12:52
% Remove the next 7 % at beginnings of lines.
\input amstex
\input amsppt.sty   
\magnification=\magstep1
\hsize 30pc
\vsize 47pc
\def\nmb#1#2{#2}         % used for renumbering, TeX should ignore.
\def\totoc{}             %= to table of content, invoked by kms-book.sty
\def\idx{}               % for producing index, invoked by kms-book.sty
\def\ign#1{}             %=ignore, invisible entry for the index only

%\pageno323
\redefine\o{\circ}

\define\ga{\gamma}

\define\ze{\zeta}

\define\th{\theta}

\define\rh{\rho}
\define\si{\sigma}
\define\ta{\tau}
\define\ph{\varphi}
\define\ch{\chi}
\define\ps{\psi}
\define\om{\omega}

\define\La{\Lambda}

\define\Ph{\Phi}
\define\Ps{\Psi}
\define\Om{\Omega}

\define\row#1#2#3{#1_{#2},\ldots,#1_{#3}}
\define\x{\times}
\define\End{\operatorname{End}}
\define\Fl{\operatorname{Fl}}

\redefine\L{{\Cal L}}
\define\ddt{\left.\tfrac \partial{\partial t}\right\vert_0}
\define\g{{\frak g}}
\def\today{\ifcase\month\or
 January\or February\or March\or April\or May\or June\or
 July\or August\or September\or October\or November\or December\fi
 \space\number\day, \number\year}
\topmatter
\title  Characteristic classes for $G$-structures
\endtitle
\author D\. V\. Alekseevsky \\
Peter W\. Michor  \endauthor
\affil
Center `Sophus Lie', Moscow \\
Institut f\"ur Mathematik, Universit\"at Wien, Austria
\endaffil
\address 	 D\. V\. Alekseevsky: Center `Sophus Lie', 
Krasnokazarmennaya 6, 111250 Moscow, USSR
\endaddress
\address
P\. W\. Michor: Institut f\"ur Mathematik, Universit\"at Wien,
Strudlhofgasse 4, A-1090 Wien, Austria
\endaddress
\email MICHOR\@AWIRAP.BITNET \endemail
\date {\today} \enddate
\thanks Supported by Project P 7724 PHY
of `Fonds zur F\"orderung der wissenschaftlichen Forschung'.
\endthanks
\abstract 
Let $G\subset GL(V)$ be a linear Lie group with Lie algebra $\g$ 
and let $A(\g)^G$ be the subalgebra of $G$-invariant elements of 
the associative supercommutative algebra $A(\g)= 
S(\g^*)\otimes \La(V^*)$. To any $G$-structure $\pi:P\to M$ 
with a connection 
$\om$ we associate a homomorphism $\mu_\om:A(\g)^G\to \Om(M)$.
The differential forms $\mu_\om(f)$ for $f\in A(\g)^G$ which are 
associated to the $G$-structure $\pi$ can be used to construct 
Lagrangians. If $\om$ has no torsion the differential forms 
$\mu_\om(f)$ are closed and define characteristic classes of a 
$G$-structure. The induced homomorphism 
$\mu'_\om:A(\g)^G\to H^*(M)$ does not depend on the choice of 
the torsionfree connection $\om$ and it is the natural generalization 
of the Chern Weil homomorphism.
\endabstract
\endtopmatter

%\input amspptb.sty
%\userunningheads
%\def\leftheadtext{\smc Peter W. Michor}
%\def\rightheadtext{\smc }
%\def\bottremark{\today\hfill}

\document

%\heading Table of contents \endheading
%\input \jobname.toc
%\loadtoc
%\loadindex

\heading\totoc\nmb0{1}. $G$-structures \endheading

\subheading{\nmb.{1.1}. $G$-structures} By a $G$-structure on a 
smooth finite dimensional manifold $M$ we mean a principal fiber 
bundle $\pi:P\to M$ together with a representation $\rh: G\to GL(V)$ 
of the structure group
in a real vector space $V$ of dimension $\dim M$ and a 1-form $\si$ 
(called the \idx{\it soldering form}) on $M$ with values in the 
associated bundle $P[V,\rh]=P\x_G V$ which is fiber wise an 
isomorphism and identifies $T_xM$ with $P[V]_x$ for each $x\in M$. 
Then $\si$ corresponds uniquely to a $G$-equivariant 1-form 
$\th\in\Om^1_{\text{hor}}(P;V)^G$ which is strongly horizontal in the 
sense that its kernel is exactly the vertical bundle $VP$. The form 
$\th$ is called the \idx{\it displacement form} of the $G$-structure.
A $G$-structure is called \idx{\it 1-integrable} if it admits 
torsionfree connections, see \nmb!{1.4} below.

We fix this setting $((P,p,M,G),(V,\rh),\th)$ from now on.

\subheading{\nmb.{1.2}. Invariant forms} We consider a multilinear form 
$f\in \bigotimes^k V^* = L^k(V)$ which is invariant in the 
sense that $f\o (\bigotimes^k\rh(g))= f$ for each $g\in G$.
Let us denote by $L^k(V)^G$ the space of all these invariant forms.
For each $f\in L^k(V)^G$ we have for any $X\in \g$, the Lie 
algebra of $G$, 
$$\align
0&=\tfrac d{dt}|_0 f(\row {\rh(\exp(tX))v}1k), \\
&= \sum_{i=1}^kf(v_1,\dots,\rh'(X)v_i,\dots,v_k),
\endalign$$
where $\rh'=T_e\rh:\g\to \g\frak l(V)$ is the differential 
of the representation $\rh$.

\subheading{\nmb.{1.3} Products of differential forms}
For $\ph\in\Om^p(P;\g)$ and $\Ps\in\Om^q(P;V)$ let us define the 
form	$\rh'_\wedge (\ph)\Ps\in\Om^{p+q}(P;V)$	by 
$$\multline
(\rh'_\wedge (\ph)\Ps)(\row X1{p+q}) =\\
= \frac 1{p!\,q!} \sum_{\si} \text{sign}(\si)
 	\rh'(\ph(\row X{\si1}{\si p}))\Ps(\row X{\si(p+1)}{\si(p+q)}).
\endmultline$$
Then $\rh'_\wedge (\ph):\Om^*(P;V)\to \Om^{*+p}(P;V)$ is a graded 
$\Om(P)$-module homomorphism of degree $p$.
Recall also that $\Om(P;\g)$ is a graded Lie algebra with the bracket
$$\multline [\ph,\ps]_{\wedge}(\row X1{p+q}) = \\
	= \frac 1{p!\,q!} \sum_{\si} \text{sign}\si\,
	[\ph(\row	X{\si1}{\si p}),\ps(\row X{\si(p+1)}{\si(p+q)})]_{\g}.
\endmultline $$
One may easily check that for the graded commutator in 
$\End(\Om(P;V))$ we have
$$\rh'_\wedge ([\ph,\ps]_\wedge ) = 
[\rh'_\wedge (\ph),\rh'_\wedge (\ps)] = 
\rh'_\wedge (\ph)\o \rh'_\wedge (\ps) - (-1)^{pq} 
\rh'_\wedge (\ps)\o \rh'_\wedge (\ph)$$
so that $\rh'_\wedge :\Om^*(P;\g) \to \End^*(\Om(P;V))$ is a 
homomorphism of graded Lie algebras.

Let $\bigotimes V$ be the tensoralgebra generated by $V$.
For $\Ph,\Ps\in \Om(P;\bigotimes V)$ we will use the associative 
bigraded product
$$\multline
(\Ph\otimes_\wedge \Ps)(\row X1{p+q}) = \\
= \frac 1{p!\,q!} \sum_{\si} \text{sign}(\si)
 	\Ph(\row X{\si1}{\si p})\otimes \Ps(\row X{\si(p+1)}{\si(p+q)})
\endmultline$$

\subheading{\nmb.{1.4}. The covariant exterior derivative} Let 
$\om\in\Om^1(P;\g)^G$ be a principal connection on the principal 
bundle $(P,p,M,G)$. Let $\ch:TP\to HP$ denote the corresponding 
projection onto the horizontal bundle $HP:=\ker \om$. The covariant 
exterior derivative $d_\om:\Om^k(P;V)\to \Om^{k+1}_{\text{hor}}(P;V)$ 
is then given as usual by 
$d_\om\Ps=\ch^*d\Ps = (d\Ps)\o \La^{k+1}(\ch)$.

\proclaim{Lemma}
For $\Ps\in \Om_{\text{hor}}(P;V)^G$ the covariant exterior derivative 
is given by $d_\om\Ps = d\Ps + \rh'_\wedge (\om)\Ps$.
\endproclaim

\demo{Proof}
If we insert one vertical vector field, say the fundamental vector 
field $\ze_X$ for 
$X\in \g$, into $d_\om\Ps$, we get 0 by definition. For the 
right hand side we use $i_{\ze_X}\Ps=0$ and 
$\L_{\ze_X}\Ps= \ddt (\Fl^{\ze_X}_t)^*\Ps = \ddt \Ps\o \La^p(r^{\exp tX}) = 
\ddt \rh(\exp(-tX))\Ps = -\rh'(X)\Ps$ to get
$$\align
i_{\ze_X}(d\Ps + \rh'_{\wedge}(\om)\Ps) &= i_{\ze_X}d\Ps + 
	di_{\ze_X}\Ps + \rh'_{\wedge}(i_{\ze_X}\om)\Ps - 
     \rh'_{\wedge}(\om)i_{\ze_X}\Ps\\
&= \L_{\ze_X}\Ps + \rh'_{\wedge}(X)\Ps = 0.
\endalign$$
Let now all vector fields $\xi_i$ be horizontal, then we get
$$\gather
(d_\om\Ps)(\row \xi 0k) = (\ch^*d\Ps)(\row \xi 0k) = 
	d\Ps(\row \xi0k),\\
(d\Ps + \rh'_{\wedge}(\om)\Ps)(\row\xi0k) = d\Ps(\row \xi0k).\qed
\endgather$$
\enddemo

\definition{\nmb.{1.5}. Definition}
If $\th\in \Om^1_{\text{hor}}(P;V)^G$ is the displacement form of a 
$G$-structure then the \idx{\it torsion} of the connection $\om$ with 
respect to the $G$-structure is $\ta:=d_\om\th= d\th+\rh'_{\wedge}(\om)\th$.
\enddefinition
Recall that a $G$-structure is called 1-integrable if it admits a 
connection without torsion. This notion has also been investigated in 
\cite{Kol\'a\v r, Vadovi\v cov\'a} where it was called prolongable.

\subheading{\nmb.{1.6}. Chern-Weil forms}
For differential forms $\ps_i\in\Om^{p_i}(P;V)$ and 
$f\in L^k(V)=(\bigotimes^k V)^*$ 
we can construct the following differential forms:
$$\gather
\ps_1\otimes_\wedge \dots \otimes_\wedge \ps_k 
	\in \Om^{p_1+\dots+p_k}(P;V\otimes \dots \otimes V),\\
f^{\ps_1,\dots,\ps_k}:= f\o (\ps_1\otimes_\wedge \dots \otimes_\wedge \ps_k ) 
	\in \Om^{p_1+\dots+p_k}(P).
\endgather$$
The exterior derivative of the latter one is clearly given by
$$\multline
d( f\o(\ps_1\otimes_\wedge\dots\otimes_\wedge\ps_k)) = 
	f\o d(\ps_1\otimes_\wedge\dots\otimes_\wedge\ps_k) =\\
= f\o \left( \tsize\sum_{i=1}^k(-1)^{p_1+\dots+p_{i-1}} 
  	\ps_1\otimes_\wedge\dots\otimes_\wedge d\ps_i\otimes_\wedge
	\dots\otimes_\wedge\ps_k \right).
\endmultline$$
We also set $f^\ps:= 
f^{\ps,\dots,\ps}=\operatorname{alt}f(\ps,\dots,\ps)$
for $\ps\in\Om^p(P;V)$.
Note that the form $f^{\ps_1,\dots,\ps_k}$ is $G$-invariant and 
horizontal if all $\ps_i\in \Om^{p_i}_{\text{hor}}(P;V)^G$ and 
$f\in L^k(V)^G$. It is then the pullback of a form on $M$.

\proclaim{\nmb.{1.7}. Lemma}	Let $0\ne \ps\in \Om^p(P;V)$ and 
$f\in L^k(V)$. Then we have:
$$f^\ps\ne 0 \iff \cases \operatorname{alt} f \ne 0, 
						&\text{ if $p$ is odd,}\\
					\operatorname{sym} f \ne 0, &\text{ if $p$ 
	                         is even,}\endcases$$
where $\operatorname{alt}$ and $\operatorname{sym}$ are the natural 
projections onto $\La(V^*)$ and $S(V^*)$, respectively. \qed
\endproclaim

\proclaim{\nmb.{1.8}. Lemma} If $f\in L^k(V)^G$ is invariant 
then we have
$$f\o \left( \tsize\sum_{i=1}^k(-1)^{p_1+\dots+p_{i-1}} 
  	\ps_1\otimes_\wedge\dots\otimes_\wedge \rh'_\wedge(\om)\ps_i
	\otimes_\wedge	\dots\otimes_\wedge\ps_k \right) = 0. $$
\endproclaim 

\demo{Proof}
This follows from the infinitesimal condition of invariance for $f$ 
given in \nmb!{1.2} by applying the alternator.
\qed\enddemo

\heading\totoc\nmb0{2}. Obstructions to 1-integrability of 
$G$-structures \endheading

\proclaim{\nmb.{2.1}. Proposition} Let $\pi:P\to M$ be a 
$G$-structure and let $f\in L^k(V)^G$ be an invariant tensor. For 
arbitrary $G$-equivariant horizontal $V$-valued forms 
$\ps_i\in\Om^{p_i}_{\text{hor}}(P;V)^G$ we consider the 
$(p_1+\dots+p_k)$-form $f^{\ps_1,\dots,\ps_k}$ on $M$ as above. 
If there is a connection $\om$ for the $G$-structure 
$\pi$ such that $d_\om\ps_i=0$ for all $i$, then the form 
$f^{\ps_1,\dots,\ps_k}$ is closed.
\endproclaim

\demo{Proof}
We use $d_\om\ps_i= d\ps_i + \rh'_\wedge(\om)\ps_i$ from lemma 
\nmb!{1.4}, and lemma \nmb!{1.8}, to obtain
$$d f^{\ps_1,\dots,\ps_k}= 
	f\o \left( \tsize\sum_{i=1}^k(-1)^{p_1+\dots+p_{i-1}} 
  	\ps_1\otimes_\wedge\dots\otimes_\wedge d_\om\ps_i
	\otimes_\wedge	\dots\otimes_\wedge\ps_k \right)= 0.\qed$$
\enddemo

\proclaim{\nmb.{2.2}. Corollary} 
1. For a $G$-structure $\pi:P\to M$ 
with displacement form $\th$ we 
have a natural homomorphism of associative algebras
$$\gather
\nu: \La(V^*)^G \to \Om(M),\\
f\mapsto f^\th=f(\th,\dots,\th).
\endgather$$
2. If the $G$-structure is 1-integrable then the image of $\nu$ consists 
of closed forms and we get an induced homomorphism
$$\nu^*:\La(V^*)^G\to H^*(M).$$
If $M$ and $G$ are compact then $\nu^*$ is injective. 
\endproclaim

\demo{Proof}
If the $G$-structure is 1-integrable then there is a connection $\om$ 
with vanishing torsion $\ta=d_\om\th=0$. Then the result follows from 
proposition \nmb!{2.1}. 

If $G$ is compact, any torsionfree connection $\om$ for $\pi:P\to M$ 
is the Levi-Civita connection for some Riemannian metric. Any form 
$f^\th$, which is parallel with respect to to $\om$, is harmonic and 
can thus not be exact for compact $M$. So $\nu^*$ is injective.
\qed\enddemo

\remark{Problem}
Is the homomorphism $\nu^*$ injective for compact $M$ but noncompact $G$?
\endremark

\subhead\nmb.{2.3}. Remark \endsubhead
Given a principal connection $\om$ on $P$ there is the induced 
covariant exterior derivative 
$\nabla:\Om^p(M;P[V])\to \Om^{p+1}(M;P[V])$ on the associated vector 
bundle $P[V]$. The soldering form (see \nmb!{1.1}) $\si:TM\to P[V]$ 
is an isomorphism of vector bundles and we may consider the pull back 
covariand derivative $\si^*\nabla$ on $TM$. Next we consider the 
`combined' covariant derivative $D^{\si^*\nabla,\nabla}$ on the 
vector bundle $L(TM,P[V])$ given by 
$D^{\si^*\nabla,\nabla}_XA=\nabla_X\o A-A\o (\si^*\nabla)_X$. 
Obviously we have $D^{\si^*\nabla,\nabla}\si=0$. Consequently 
for any $f\in L^k(V)^G$ we have that $f^\th \in \Om^k(M)$ is parallel 
for the connection induced on $\La^kT^*M$ from $\si^*\nabla$ on $TM$.

\heading\totoc\nmb0{3}. The generalized Chern-Weil homomorphism for 
$G$-structures \endheading

\subheading{\nmb.{3.1}. The Chern-Weil homomorphism}
Let $\om$ be a connection for a $G$-structure $\pi:P\to M$ with 
curvature form $\Om\in\Om^2_{\text{hor}}(P,\g)$. Then the 
Bianchi identity $d_\om\Om=0$ holds. If we apply proposition 
\nmb!{2.1} to $\ps_i=\Om$ we obtain a homomorphism
$$\ga: S(\g^*)^G\to \Om(M),$$ 
given by $\ga(f)=f^\Om$. Since the image of $\ga$ consists of closed 
forms we have an induced homomorphism 
$$\ga':S(\g^*)^G\to H^*(M).$$
This is the well known Chern-Weil homomorphism. 

\subheading{\nmb.{3.2}. The algebra $A(\g,V)$}
In oder to generalize the Chern Weil homomorphism we associate to a 
Lie algebra $\g$ and a vector space $V$ the 
associative graded commutative algebra 
$$A(\g,V):=S(\g^*)\otimes \La(V^*),$$
where the generators of the symmetric algebra $S(\g^*)$ have degree 2.
We may also consider $A(\g,V)$ as a graded Lie algebra with the 
bracket 
$$[a\otimes \ph,b\otimes \ps]:=\{a,b\}\otimes \ph\wedge \ps,
	\quad a,b\in S(\frak g^*),\ph,\ps\in\La(V^*),$$
where $\{a,b\}$ is the usual Poisson-Lie bracket in $S(\g^*)$.

Let now $\g$ be the Lie algebra of the Lie group $G$ and let 
$\rh:G\to GL(V)$ be a representation. 
Then $G$ acts naturally on $A(\g,V)$, and we denote $A(\g,V)^G$ the 
subalgebra of $G$-invariant elements in $A(\g,V)$.

\subheading{\nmb.{3.3}. Remark} The associative algebra $A(\g,V)^G$ 
contains the subalgebra $S(\g^*)^G\otimes \La(V^*)^G$, in general as 
a proper subalgebra. Actually, let $G\subset GL(V)$ be the isotropy 
group of an irreducible Riemannian symmetric space $M$. Then the 
curvature tensor of $M$ defines an element of 
$(\g^*\otimes\La^2V^*)^G\subset A(\g,V)^G$ that does not belong to 
$(\g^*)^G\otimes (\La^2V^*)^G=0$

\subheading{\nmb.{3.4}. The generalized Chern-Weil homomorphism}
Now we are in a position to combine the constructions \nmb!{2.2} and 
\nmb!{3.1}.

\proclaim{Theorem} Let $\pi:P\to M$ be a $G$-structure on $M$ with 
displacement form $\th$. Any connection $\om$ in $\pi$ defines a 
homomorphism of associative algebras
$$\gather
\mu:A(\g,V)^G\to \Om(M)\\
(S^p(\g^*)\otimes\La^q V^*)^G\ni f\mapsto f^{\Om,\th}
=f(\undersetbrace p \to{\Om,\dots,\Om},\undersetbrace q 
\to {\th,\dots,\th})
\endgather$$
If the connection $\om$ has no torsion then the image of $\mu$ 
consists of closed forms and $\mu$ induces a homomorphism
$$ \mu':A(\g,V)^G\to H^*(M),$$
which is independent of the choice of the torsionfree connection.
\endproclaim

In other words, any $G$-invariant tensor 
$f\in S^p(\g^*)\otimes\La^q(V^*)$ defines a cohomology class 
$[f^{\Om,\th}]\in H^{2p+q}(M)$ which is an invariant of the 
1-integrable $G$-structure. We call it a 
\idx{\it characteristic class} of the 1-integrable $G$-structure $\pi$.

\demo{Proof}
It just remains to show that the cohomology class $[f^{\Om,\th}]$ 
does not depend on the choice of the torsionfree connection for the 
$G$-structure $\pi:P\to M$. 

So let $\om_0$, $\om_1$ be two torsionfree connections for the 
$G$-structure, let $\ph=\om_1-\om_0$, and denote by 
$\Om_t=d_{\om_t}\Om_t$ the curvature form of the torsionfree 
connection $\om_t=\om_0+t\ph=(1-t)\om_0+t\om_1$.
We claim that for $f\in(S^p(\g^*)\otimes\La^q(V^*))^G$ we have
$$\gather
f^{\Om_1,\th}-f^{\Om_0,\th}= d(Tf),\quad\text{ where }\tag1\\
Tf = p \int_0^1 f(\ph,\Om_t,\dots,\Om_t,\th,\dots,\th)\,dt
\endgather$$
is the transgression form of $f$ on $P$.
The assertion is immediate from \thetag1.
To prove it we compute $\partial_t f^{\Om_t,\th}$ using the 
identities $\partial_t\Om_t=d_{\om_t}\ph$ (see \cite{Kobayashi, 
Nomizu II, p\. 296}), $d_{\om_t}\Om_t=0$, and $d_{\om_t}\th=0$.
$$\align
\partial_t f^{\Om_t,\th} 
	&= p\,f(\partial_t\Om_t,\Om_t,\dots,\Om_t,\th,\dots,\th) \\
	&= p\,f(d_{\om_t}\ph,\Om_t,\dots,\Om_t,\th,\dots,\th) \\
	&= p\,d_{\om_t}f(\ph,\Om_t,\dots,\Om_t,\th,\dots,\th) \\
	&= p\,d\,f(\ph,\Om_t,\dots,\Om_t,\th,\dots,\th).\qed
\endalign$$
\enddemo

\subheading{\nmb.{3.5}. Remarks about secondary characteristic 
classes}
If the characteristic forms $f^{\Om_1,\th}$ and $f^{\Om_0,\th}$ 
associated with two torsionfree connections $\om_1$ and $\om_0$ 
vanish we obtain a secondary characteristic class $[Tf]$. It is a 
natural generalization of the classical Chern-Simons characteristic 
class, see \cite{Chern, Simons}, \cite{Kobayashi, Ochiai}.

Problem: study conditions when the secondary characteristic class 
$[Tf]$ does not depend on the choice of the torsionfree connections 
$\om_1$ and $\om_0$.

\subheading{\nmb.{3.6}. Examples of characteristic classes}
Assume that a linear group $G\subset GL(V)$ preserves some pseudo 
Euclidean metric in $V=\Bbb R^n$. Then we may identify the Lie 
algebra $\g=\operatorname{Lie}(G)$ with a subspace $\g\subset\La^2V$. 
Suppose that there exists a $G$-invariant supplement $\frak d$ to 
$\g$ in $\La^2V$. Then the $G$-equivariant projection $\La^2V\to\g$ 
along $\frak d$ determines a $G$-invariant element 
$q\in\g\otimes\La^2V^*\cong \g^*\otimes\La^2V^*.$
The element $q$ defines a 4-form $q^{\Om,\th}$ on the base of any 
$G$-structure $\pi:P\to M$ with a connection $\om$ and curvature 
$\Om$. It may be written as
$$q^{\Om,\th}=q(\Om,\th,\th)
	=q^a_{bcd}R^b_{aef}\th^c\wedge \th^d\wedge \th^e\wedge \th^f,$$
where $(q^a_{bcd})$ is the coordinate expression of $q$ in the 
standard basis $(e_a)$ of $V=\Bbb R^n$, $\th=e_a\otimes\th^a$, and 
$\Om=R^a_{bef}\th^e\wedge \th^f$.

If $\om$ is torsionfree the 4-form $q^{\Om,\th}$ is closed and it 
defines a cohomology class $[q^{\Om,\th}]\in H^4(M)$ independently of 
the choice of $\om$.

\subheading{\nmb.{3.7}. Remarks about the classification of 
characteristic classes}
The classification of characterictic classes for $G-structures$ with 
a given Lie group $G$ reduces to the construction of generators of 
the associative algebra $A(\g,V)^G=(S(\g^*)\otimes\La(V^*))^G$.
We may also use the bracket to multiply characteristic classes.
It suffices to solve this problem for those Lie groups $G$ which 
appear as holonomy groups of torsionfree connection. Only for such 
groups $G$ there exist 1-integrable non-flat $G$-structures. Under 
the additional hypothesis of irreducibility, all such groups were 
classified by \cite{Berger}, up to some gaps which were filled by 
\cite{Bryant} and \cite{Alekseevsky, Graev}.

\enddocument